\documentclass[a4paper,12pt]{article}
\usepackage[latin1]{inputenc}
\usepackage[T1]{fontenc}
\usepackage[english]{babel}
\usepackage{amssymb}
\usepackage{latexsym}

\usepackage{amsmath}
\usepackage{amsthm}
\usepackage{amsxtra}
\usepackage{eucal}
\providecommand{\LyX}{L\kern-.1667em\lower.25em\hbox{Y}\kern-.125emX\@}

\def\R{\mathbb{R}}
\def\E{\mathbb{E}}
\def\dint{\displaystyle\int}

\newtheorem{thm}{Theorem}[section]

\newtheorem{proposition}{Proposition}[section]

\newtheorem{remark}{Remark}[section]

\newtheorem{defi}[thm]{Definition}

\newtheorem{lem}[thm]{Lemma}

\newtheorem{prop}[thm]{Proposition}

\newenvironment{class}[1][AMS Classification]{\textbf{#1.} }{}
\newenvironment{MC}[1][Key words]{\textbf{#1.} }{}

\def\disp{\displaystyle{}}

\def\al{\alpha}
\def\si{\sigma}
\def\be{\beta}
\def\om{\omega}
\newcommand{\ig}[2]{\displaystyle \int_{#1}^{#2}}

\title{\textbf{On one-dimensional stochastic differential equations involving
the maximum process}}

\author{ R. BELFADLI $^{1}\quad$ S. HAMAD\`ENE   $^{2}\quad$ Y. OUKNINE $^{1}$ $^{3}\quad$
\vspace*{0.1in}\\
$^{1}$  Department of Mathematics, Faculty
of Sciences Semlalia,\\ Cadi Ayyad University 2390 Marrakesh, Morocco.\\
r.belfadli@ucam.ac.ma \quad ouknine@ucam.ac.ma
\vspace*{0.1in}\\
$^{2}$ Laboratoire de Statistique et Processus,\\Université du Mans, 72085 Le Mans Cedex 9, France.\\
hamadene@univ-lemans.fr
\vspace*{0.1in}\\
$^{3}$ Hassan II Academy of sciences and Technology\\ Rabat, Morocco.\\
}

\begin{document}
\maketitle \textbf{Abstract.} We prove existence and pathwise uniqueness results for four different types of
 stochastic differential equations (SDEs) perturbed by the past maximum process and/or the local time at zero. Along the first three studies, the coefficients are no longer Lipschitz. The first type is the equation
\begin{equation*}
\label{eq1} X_{t}=\int_{0}^{t}\sigma
(s,X_{s})dW_{s}+\int_{0}^{t}b(s,X_{s})ds+\alpha \max_{0\leq s\leq
t}X_{s}.
\end{equation*}
The second type is the equation
\begin{equation*}
\label{eq2} \left\{
\begin{array}{l}
X_{t} =\ig{0}{t}\sigma (s,X_{s})dW_{s}+\ig{0}{t}b(s,X_{s})ds+\alpha
\max_{0\leq s\leq
t}X_{s}\,\,+L_{t}^{0},\,\,\, \\
X_{t} \geq 0,\,\forall t\geq 0.
\end{array}
\right.
\end{equation*}
The third type is the equation
\begin{equation*}
\label{eq3} X_{t}=x+W_{t}+\int_{0}^{t}b(X_{s},\max_{0\leq u\leq
s}X_{u})ds.
\end{equation*}
We end the paper by establishing the existence of strong solution and pathwise uniqueness, under Lipschitz condition, for the SDE
\begin{eqnarray*}\label{e2}
X_t=\xi+\int_0^t \si(s,X_s)dW_s +\int_0^t b(s,X_s)ds +\al\max_{0\leq
s\leq t}X_s +\be \min_{0\leq s \leq t}X_s.
\end{eqnarray*}

\begin{MC}
Perturbed stochastic differential equations; Strong solution; Pathwise uniqueness; Local time.
\end{MC}

\begin{class}
Primary 60H10; Secondary 60J60.
\end{class}
\section{Introduction}
One-dimensional stochastic differential equations (SDEs) involving
the past maximum process and/or the past minimum process and/or the local time process has attracted several
authors (see for instance \cite{dz},  \cite{davis1}, \cite{ChDo}, \cite{ChDoHu}, \cite{PrWe}).

In this paper, we will study four classes of these SDEs, three with singular coefficients and the other with Lipschitz one. In fact, we are concerned with the existence and pathwise uniqueness of solutions. Because of the lack of regularity of those latter, the usual fixed point approach cannot be applied.
So our method is based on uniqueness in law and basic results about local times,
namely the Tanaka's formula. We do not require strong regularity assumption on the drift coefficients. Besides, we allow  the diffusion coefficient be discontinuous if it is strictly positive.

Let $\Omega$ be the set of continuous functions from $R^+$ into $R$,
$P$ the Wiener measure on $\Omega$, $(W_t)_{t\geq 0}$ the process of
coordinate maps from $\Omega$ into $R$, ${\cal
F}=\sigma\{W_t,\,t\geq 0\}$, $({\cal F}_t)_{t\geq 0}$ the completion
of the natural filtration of $W$ with the $P$-null sets of $\cal F$.
Therefore $(W_t)_{\geq 0}$ is a standard Brownian motion on the
filtered probability space $(\Omega, {\cal F}, ({\cal F}_t)_{t\geq
0},P)$.

In Section \ref{section1}, we establish both the existence of a strong solution and the pathwise uniqueness for the SDE
\begin{equation}
\label{eq1} X_{t}=\int_{0}^{t}\sigma
(s,X_{s})dW_{s}+\int_{0}^{t}b(s,X_{s})ds+\alpha \max_{0\leq s\leq
t}X_{s}.
\end{equation}

In Section \ref{section2}, we study the perturbed SDEs with reflecting boundary
\begin{equation}
\label{eq2} \left\{
\begin{array}{l}
X_{t} =\ig{0}{t}\sigma (s,X_{s})dW_{s}+\ig{0}{t}b(s,X_{s})ds+\alpha
\max_{0\leq s\leq
t}X_{s}\,\,+L_{t}^{0},\,\,\, \\
X_{t} \geq 0,\,\forall t\geq 0.
\end{array}
\right.
\end{equation}
 Actually to equation (\ref{eq1}) we add a reflecting process which is the local time
$(L_{t}^{0})_{t\geq 0}$ of the process $X$ at $0$. Its role is to
push upward the process $X$ in order to keep it above $0$, i.e. to
have the condition $X\geq 0$ satisfied.\\
Section \ref{section3} is devoted to the study of SDEs with a drift containing the maximum process
\begin{equation}\label{eq3}
X_{t}=x+W_{t}+\int_{0}^{t}b(X_{s},\max_{0\leq u\leq s}X_{u})ds.
\end{equation}\\
Once more we show existence and uniqueness of the solution
$(X_t)_{t\geq 0}$. Note that $\left( X_{t}\right) _{t\geq 0}$ is not
a Markov process but $(X_{t},\max_{0\leq s\leq t}X_{s})_{t\geq 0}$
is a Markov process.

In Section \ref{section4} we deal with the doubly perturbed SDEs
\begin{eqnarray}\label{e22}
X_t=\xi+\int_0^t \si(s,X_s)dW_s +\int_0^t b(s,X_s)ds +\al\max_{0\leq
s\leq t}X_s +\be \min_{0\leq s \leq t}X_s.
\end{eqnarray}
Under Lipschitz condition on the coefficients and an assumption on the parameters $\al$, $\be$, we prove the strong existence and pathwise uniqueness property for (\ref{e22}). We conclude this paper by a discussion
on some further extensions of our results. $\diamond$
\bigskip

Throughout this paper, solutions to each stochastic differential
equation under consideration should be understood as continuous
processes. Besides let us recall the following definitions:
\medskip
A {\it strong solution} for (\ref{eq1}) (resp. (\ref{eq2}); resp.
(\ref{eq3})) on the probability space $(\Omega, \mathcal{F}, P)$
endowed with the Brownian motion $W$ is a continuous process
$(X_t)_{t\geq 0}$ adapted w.r.t. the natural filtration of $W$ and
which satisfies (\ref{eq1}) (resp. (\ref{eq2}); resp. (\ref{eq3})).
A {\it weak solution} for (\ref{eq1}) (resp. (\ref{eq2}); resp.
(\ref{eq3})) on a filtered probability space
$(\Omega,\mathcal{F},(\mathcal{F}_{t}) _{t\geq 0},P)$ is a couple of
processes $(X,W)$ on that space such that $X$ is adapted with
respect to $(\mathcal{F}_{t})_{t\geq 0}$, $W$ is a Brownian motion
adapted to $(\mathcal{F}_{t}) _{t\geq 0}$ and finally the couple
$(X,W)$ satisfies (\ref{eq1}) (resp. (\ref{eq2}); resp.
(\ref{eq3})). We say that {\it uniqueness in law} holds for
(\ref{eq1}) (resp. (\ref{eq2}); resp. (\ref{eq3})), if any two weak
solutions $(X,W)$ and $(\tilde{X},\tilde{W})$ have the same laws
whenever the laws at the initial time $t=0$ are the same.

Finally we say that {\it pathwise uniqueness} holds for equation
(\ref{eq1}) (resp. (\ref{eq2}); resp. (\ref{eq3})) if two weak
solutions $(X,W)$ and $(\tilde{X},W)$ defined on the same filtered
probability space with the same Brownian motion $W$ are such that
the processes $X$ and $\tilde{X}$ being undistinguishable.
$\diamond$
\medskip

There were many works which discuss under which conditions on the
drift $b$ and the diffusion coefficient $\sigma$ we have the
existence of strong solutions of stochastic differential equations.
In the case when the equation is one-dimensional and $\sigma$ is not
degenerated, several results have been obtained by Y.Ouknine
(\cite{ou1},\cite{ou2},\cite{ou3}). For SDEs which involves the local
time of the unknown process $X$, as far as we know, the most general
result is given by M.Rutkowski in $(\cite{ru})$ where he showed that
the so called $\mathbf{LT}$-condition is sufficient to have pathwise
uniqueness. In that case the diffusion coefficient $\sigma$ can be
degenerated. $\diamond$
\section{Perturbed SDEs with measurable coefficients}\label{section1}
To begin with, let us recall the following definition which we
borrow from Ash $\&$ Dol\'{e}ans-Dade $\cite{ad}$.
\begin{defi}
A two-variable function $f(t,x)$ is called monotonically increasing if whenever
$t_{2}\geq t_{1}$ and $x_{2}\geq x_{1}$ we have:
\begin{equation*}
f(t_{2},x_{2})-f(t_{2},x_{1})-f(t_{1},x_{2})+f(t_{1},x_{1})\geq 0.
\end{equation*}
\end{defi}
To a monotonically increasing function $f\left( t,x\right) $ we  associate the Lebesgue-Stieltjes measure defined, for any $t_{2}\geq t_{1}$ and $x_{2}\geq x_{1}$, by:
\begin{equation*}
\mu \left( \lbrack t_{1},t_{2})\times \lbrack x_{1},x_{2})\right)
=f(t_{2},x_{2})-f(t_{2},x_{1})-f(t_{1},x_{2})+f(t_{1},x_{1}).
\end{equation*}
Thus, for a measurable function $g\left( t,x\right)$, we can define its
Lebesgue-Stieltjes integral w.r.t. $f$ in the following way:
\begin{equation*}
\int_{t_{1}}^{t_{2}}\int_{a}^{b}g\left( t,x\right)
d_{t,x}f(t,x)=\int_{t_{1}}^{t_{2}}\int_{a}^{b}g\left( t,x\right)
\mu(dt,dx).
\end{equation*}
Let us now make precise the assumptions on the coefficients $b$ and
$\sigma$: \medskip

$(i)$ $\sigma$ and $b$ are measurable functions on
$\mathbf{R}^{+}\mathbf{\times R\,}$ and sub-linearly growing, i.e.,
\begin{equation}
\left| \sigma (t,x)\right| +\left| b(t,x)\right| \leq C\left(
1+\left| x\right| \right), \forall t\in \mathbf{R}^{+}\mbox{ and
}x\in \mathbf{R}. \tag{LG}
\end{equation}

$(ii)$ There exists a monotonically increasing function $f\left(
t,x\right) $ on $ \mathbf{R}^{+}\times \mathbf{R}$ such that
$f\left( 0,x\right) $ is increasing and
\begin{equation}
\left| \sigma (t,x)-\sigma (t,y)\right| ^{2}\leq \left| f\left(
t,x\right) -f\left( t,y\right) \right| ,\forall t\in \mathbf{R} ^{+}
\mbox{ and }x,y\in \mathbf{R}.  \tag{BV2}
\end{equation}
$(iii)$ There exists $\varepsilon >0$ such that:
\begin{equation}
\sigma (t,x)\geq \varepsilon \text{  for any }t\in \mathbf{R}^{+}
\mbox{ and }x\in \mathbf{R}.  \tag{ND}
\end{equation}

\begin{remark} We can easily see that the mapping $x\mapsto f\left( t,x\right) $ is increasing for any fixed $t\in
\mathbf{R}^{+}.$
\end{remark}

The main result of this section is the following
\begin{thm} Let $\alpha \in (0,1)$. Then, under the conditions (LG), (BV2) and (ND) the perturbed SDE (\ref{eq1})
has a strong solution which is moreover pathwise unique.
\end{thm}
\noindent {\bf Proof}: It will be obtained after the following three
steps.
\medskip

\noindent {\it Step 1}: {\bf Existence and uniqueness in law}
\medskip

\noindent We start with a Brownian motion $W$, and consider the process $X_{t}=W_{t}+\frac{\alpha }{1-\alpha }\max_{0\leq s\leq t}W_{s}$.
Then, $X$ solves the SDE \[X_{t}=W_{t}+\alpha \max_{0\leq s\leq t}X_{s}.\]
Therefore, using the condition (LG) and the Girsanov's theorem, we easily deduce the existence and uniqueness in law of the solution of the SDE
\begin{equation} \label{sigma1}
X_{t}=W_{t}+\int_{0}^{t}b(s,X_{s})ds+\alpha
\max_{0\leq s\leq t}X_{s}.
\end{equation}
\medskip

We now claim that uniqueness in law holds for solutions of (\ref{eq1}).
Actually since $\sigma $ is bounded from below by a real positive
constant, there is a one to one correspondence between the
distributions of solutions of (\ref{eq1}) and the ones of the
same equation but with $\sigma \equiv 1$. Indeed, if $Y$ is a
solution to (\ref{eq1}) and if we define $T_t$, for any $t\geq 0$, by
\begin{equation*}
T_{t}=\inf \left\{ s>0:\int_{0}^{s}\sigma (u,Y_{u})^{2}du\geq
t\right\}
\end{equation*}
then the process $(Y_{T_{t}})_{t\geq 0}$ is a solution to
(\ref{sigma1}), and vice versa.
\medskip

\noindent {\it Step 2}: {\bf Pathwise uniqueness}
\medskip

Let us first show that $L_{.}^{0}(X-Y)\equiv0$, whenever $X$ and $Y$ denote any two solutions of the SDE (\ref{eq1}) with
the same underlying Brownian motion $W$.
By the right continuity of $L_{.}^{0}$ it is enough to prove that, for any $t\geq 0$,
\begin{equation*}
\int_{0}^{+\infty }\frac{L_{t}^{a}\left( X-Y\right)}{a}da<+\infty.
\end{equation*}
Indeed, using the density occupation formula we can write for any $\delta >0$,
\begin{equation*}
\int_{\delta }^{+\infty }\frac{L_{t}^{a}}{a}da=\int_{0}^{t}\frac{d\langle
X-Y\rangle_{s}}{X_{s}-Y_{s}}1_{\left\{ X_{s}-Y_{s}>\delta \right\}
}=\int_{0}^{t}\frac{\left( \sigma (s,X_{s})-\sigma (s,Y_{s})\right)
^{2}}{%
X_{s}-Y_{s}}1_{\left\{ X_{s}-Y_{s}>\delta \right\} }ds.
\end{equation*}
Applying the assumption (BV2) we obtain
\begin{equation*}
\int_{0}^{t}\frac{\left( \sigma (s,X_{s})-\sigma (s,Y_{s})\right)
^{2}}{%
X_{s}-Y_{s}}1_{\left\{ X_{s}-Y_{s}>\delta \right\} }ds\leq
\int_{0}^{t}\frac{%
|f(s,X_{s})-f(s,Y_{s})|}{X_{s}-Y_{s}}1_{\left\{ X_{s}-Y_{s}>\delta
\right\} }ds.
\end{equation*}
As a consequence,
\begin{eqnarray}
  \E\left[\int_{\delta }^{+\infty }\frac{L_{t}^{a}(X-Y)}{a}da\right]\leq& \E\left[\dint_{0}^{t}\frac{
|f(s,X_{s})-f(s,Y_{s})|}{X_{s}-Y_{s}}1_{\left\{ X_{s}-Y_{s}>\delta
\right\} }ds\right].
\end{eqnarray}
Now, by a localization argument we may assume that $f$ is bounded and we consider the sequence of functions $f_n$ defined by
\begin{equation*} f_{n}(t,a)=( f(t,.)\ast \theta _{n})(a)
\end{equation*}
where $\theta _{n}$ is the standard positive regularizing mollifiers sequence. So, for fixed $t\geq 0$,
\[f_{n}(t,a)\rightarrow f(t,a) \quad \, \text{for every}\,\, a \notin D_{t},\]
where $D_{t}$ is the denumerable set of discontinuous points of the function $a\mapsto f(t,a)$.

Hence, using successively Fatou's Lemma, the intermediate value theorem, the fact that $\sigma \geq
\epsilon$ and \[\frac{d}{dt}\langle \alpha X+(1-\alpha)Y \rangle _t=(\alpha \sigma
(t,X_t)+(1-\alpha)\sigma (t,Y_t))^2 \geq \epsilon^2\] we obtain
\begin{align*}
& \E\left[\dint_{0}^{t}\dfrac{\left( f(s,X_{s})-f(s,Y_{s})\right)
}{X_{s}-Y_{s}} 1_{\left\{ X_{s}-Y_{s}>\delta \right\} }ds\right]\\
&\leq \liminf_{n\rightarrow +\infty }\E\left[\int_0^t\dfrac{\left(
f_{n}(s,X_{s})-f_{n}(s,Y_{s})\right) }{X_{s}-Y_{s}}1_{\left\{
X_{s}-Y_{s}>\delta \right\} }ds\right]\\
&= \liminf_{n\rightarrow +\infty
}\E\left[\int_{0}^{t}ds\int_{0}^{1}\dfrac{%
\partial f_{n}}{\partial a}(s,\alpha X_{s}+(1-\alpha )Y_{s})d\alpha \right]\\
&= \liminf_{n\rightarrow +\infty
}\int_{0}^{1}d\alpha\E\left[\int_{0}^{t}\dfrac{
\partial f_{n}}{\partial a}(s,\alpha X_{s}+(1-\alpha )Y_{s})ds\right]\\
&\leq \frac{1}{\epsilon^2}\liminf_{n\rightarrow +\infty
}\int_{0}^{1}d\alpha\int_R da \E\left[\int_0^t\frac{\partial
f_n}{\partial a}(s,a)dL_{s}^{a}(\alpha X+(1-\alpha )Y)\right].
\end{align*}
Note that we have used in the first inequality the fact that
\begin{eqnarray}\label{statement}
  \dint_{0}^{t}P[ (X_{s}\in D_{s})\cup (Y_{s}\in D_{s})]ds &=& 0.
\end{eqnarray}
To see that this statement holds, it suffices to consider the case where $\si\equiv1$ and $b\equiv0$, the general situation may be deduced by applying both  Girsanov theorem and a time change. In this case, as it was observed in the beginning of the proof, the solution may be expressed as a function of $(W_t,M^{W}_t)$ , and since the law of the later process is absolutely continuous w.r.t the Lebesgue measure and $D_{s}$ is denumerable we get (\ref{statement}).

Next, since $f$ is
monotonically increasing, the positivity of the measure $\dfrac{\partial
^{2}f_{n}}{\partial t\partial a}(dt,da)$ yields
\begin{align*}
&  \E\left[\int_R da\int_{0}^{1}d\alpha\int_{0}^{t}\dfrac{\partial
f_{n}}{\partial a}(s,a)dL_{s}^{a}(\alpha X+(1-\alpha )Y)\right] \nonumber\\
 & = \E\left[\int_R \int_{0}^{1}\dfrac{\partial f_{n}}{\partial
a}(t,a)L_{t}^{a}(\alpha X+(1-\alpha )Y)d\alpha da\right.\\
& \left. \quad \quad -\int_R \int_{0}^{1}d\alpha\int_{0}^{t}\dfrac{\partial ^{2}f_{n}}{\partial
t\partial a}(ds,da)L_{s}^{a}(\alpha X+(1-\alpha )Y)\right]\nonumber\\
&\leq  \E\left[\int_R \int_{0}^{1}\dfrac{\partial f_{n}}{\partial
a}(t,a)L_{t}^{a}(\alpha X+(1-\alpha )Y) d\alpha da\right].
\end{align*}
Thus,
\begin{eqnarray}
  \E\left[\int_{\delta }^{+\infty }\frac{L_{t}^{a}(X-Y)}{a}da\right] &\leq& \liminf_{n\rightarrow +\infty}\E\left[\int_R \int_{0}^{1}\dfrac{\partial f_{n}}{\partial
a}(t,a)L_{t}^{a}(\alpha X+(1-\alpha )Y) d\alpha da\right].
\end{eqnarray}

However, since $\alpha \in (0,1)$, then standard
calculations imply that, for any $p\geq 0$ and $t\geq 0$,
\begin{equation}
\label{esti1} \E[\sup_{s\leq t}|X_s|^p]<\infty.
\end{equation}
So, if $X^\alpha$ is the process defined by $X^\alpha :=\alpha X+(1-\alpha )Y$ we deduce by using Tanaka formula and the inequality
$|X_t^\alpha-a|-|X_0^\alpha-a|\leq
|X_t^\alpha-X_0^\alpha|$ that
$$
\sup_{\alpha \in [0,1],a\in R}\E \left[L_{t}^{a} (\alpha X+(1-\alpha
)Y)\right]<\infty.$$ Therefore, we obtain

\begin{eqnarray}\label{esti}\nonumber\\
  \E\left[\int_{\delta }^{+\infty }\frac{L_{t}^{a}}{a}da\right] &\leq& \sup_{\alpha
\in [0,1],a\in R}\E\left[ L_{t}^{a}
(\alpha X+(1-\alpha )Y)\right]\int_R \dfrac{\partial f_{n}%
}{\partial a}(t,a)da  \leq  C\left\| f\right\|_{\infty }
\end{eqnarray}
where \ $C>0$\ is a generic constant. 
The result follows by letting $\delta$ to zero in (\ref{esti}).

\noindent {\it Step 3}: {\bf Extrema of two solutions are solutions }
\medskip

We now show that $X\wedge Y$ and $X\vee Y$ are also solutions to
$\left( \ref {eq1}\right) $. Using Tanaka's formula we write
\begin{equation}
\label{xsupy}
\begin{array}{lll}
X_{t}\vee Y_{t} & = & Y_{t}+\left( X_{t}-Y_{t}\right) ^{+}\text{ } \\
& = & Y_{t}+\left( \int_{0}^{t}1_{\left\{ X_{s}>Y_{s}\right\}
}d\left(
X_{s}-Y_{s}\right) +\frac{1}{2}L_{t}^{0}\left( X-Y\right) \right) \\
& = & \int_{0}^{t}1_{\left\{ X_{s}>Y_{s}\right\}
}dX_{s}+\int_{0}^{t}1_{\left\{ X_{s}\leq Y_{s}\right\} }dY_{s} \\
& = & \int_{0}^{t}\sigma (s,X_{s}\vee
Y_{s})dW_{s}+\int_{0}^{t}b(s,X_{s}\vee Y_{s})ds \\
&  & \qquad \qquad +\,\alpha \left( \int_{0}^{t}1_{\left\{
X_{s}>Y_{s}\right\} }dM_{s}^{X}+\int_{0}^{t}1_{\left\{ X_{s}\leq
Y_{s}\right\} }dM_{s}^{Y}\right),
\end{array}
\end{equation}
where $M_{s}^{Z}$ denotes $\max_{0\leq s\leq t}Z_{s}$.

Next we claim that
\begin{equation*}
M_{t}^{X\vee Y}=\int_{0}^{t}1_{\left\{ X_{s}>Y_{s}\right\}
}dM_{s}^{X}+\int_{0}^{t}1_{\left\{ X_{s}\leq Y_{s}\right\}
}dM_{s}^{Y}, \forall t\geq 0.
\end{equation*}
Actually by commuting $X$ and $Y$ in (\ref{xsupy}) we obtain:
\begin{equation*} \int_{0}^{t}1_{\left\{ X_{s}=Y_{s}\right\}
}dM_{s}^{X}=\int_{0}^{t}1_{\left\{ X_{s}=Y_{s}\right\}
}dM_{s}^{Y},\, \forall t\geq 0,
\end{equation*}
and then
\begin{eqnarray*}
M_{t}^{X\vee Y} &=&\int_{0}^{t}1_{\left\{ X_{s}>Y_{s}\right\}
}dM_{s}^{X\vee
Y}+\int_{0}^{t}1_{\left\{ X_{s}\leq Y_{s}\right\} }dM_{s}^{X\vee Y} \\
&=&\int_{0}^{t}1_{\left\{ X_{s}>Y_{s}\right\}
}dM_{s}^{X}+\int_{0}^{t}1_{\left\{ X_{s}<Y_{s}\right\}
}dM_{s}^{Y}+\int_{0}^{t}1_{\left\{ X_{s}=Y_{s}\right\}
}dM_{s}^{X\vee Y}.
\end{eqnarray*}
But
\begin{equation*}
\begin{array}{lll}
\int_{0}^{t}1_{\left\{ X_{s}=Y_{s}\right\} }d( M_{s}^{X\vee
Y}-M_{s}^{Y}) & = & \int_{0}^{t}1_{\left\{ X_{s}=Y_{s}\right\}
}d( M_{s}^{X}-M_{s}^{Y}) ^{+} \\
& = & \int_{0}^{t}1_{\left\{ M_{s}^{X}>M_{s}^{Y}\right\} }1_{\left\{
X_{s}=Y_{s}\right\} }d( M_{s}^{X}-M_{s}^{Y}) =0.\\
&  &
\end{array}
\end{equation*}
Thus
\begin{equation*}
M_{t}^{X\vee Y}=\int_{0}^{t}1_{\left\{ X_{s}>Y_{s}\right\}
}dM_{s}^{X}+\int_{0}^{t}1_{\left\{ X_{s}\leq Y_{s}\right\}
}dM_{s}^{Y}, \forall t\geq 0.
\end{equation*}
This shows that $X\vee Y$ is actually a solution for (\ref{eq1}).
Similarly one can show that $X\wedge Y\,$is a solution as well. Now,
as $X$ and $Y$ have integrable paths on finite intervals and taking
into account the uniqueness in law we obtain:
\begin{equation*}
\E [|X-Y|] =\E [ X\vee Y-X\wedge Y]=0,
\end{equation*}
therefore $X=Y$.
\bigskip

By the same lines as the previous proof we have the following more
general setting:

\begin{proposition}
Assume the existence and uniqueness in law for $(\ref{eq1})$ holds
and assume also that for any solutions $X$ and $Y$ defined on the
same stochastic basis with respect to the same Brownian motion one
has $L_.^{0}\left( X-Y\right) =0$, then $(\ref{eq1})$ has a strong
solution which is pathwise unique.
\end{proposition}

\section{Perturbed SDEs with reflecting boundary}\label{section2}

In this section we focus on the stochastic differential equation
(\ref{eq2}), i.e.,
$$\left\{
\begin{array}{lll}
X_{t} &=&\int_{0}^{t}\sigma
(s,X_{s})dW_{s}+\int_{0}^{t}b(s,X_{s})ds+\alpha \max_{0\leq s\leq
t}X_{s}\,\,+L_{t}^{0},\,\,\, \\
X_{t} &\geq &0,\,\,\forall t\geq 0,
\end{array}\right.
$$
where $(L_{t}^{0})_{t\geq 0}$ is the local time of the process $X$
in $0$. To begin with, let us consider a particular case of the
previous equation, namely,

\begin{equation}
\left\{ \label{eqint1}
\begin{array}{lll}
X_{t} &=&W_{t}+\alpha \disp{\max_{0\leq s\leq
t}X_{s}}\,\,+L_{t}^{0},\\
X_{t} &\geq &0,\,\,\forall t\geq 0.
\end{array}
\right.
\end{equation}

\begin{lem}
There exists a unique solution in law for equation (\ref{eqint1}).
\end{lem}

\noindent {\bf Proof}: If we set $Y_{t}=( X_{t})^{2}$ then $Y_{t}$
satisfies the following SDE:

\begin{equation*}
Y_{t}=2\int_{0}^{t}\sqrt{Y_{s}}dW_{s}+t+\alpha \max_{0\leq s\leq
t}Y_{s}.
\end{equation*}
Actually this is due to the fact that for any $t\geq 0$,
$\{\frac{Y_t}{\max_{0\leq s\leq
t}Y_s}\}^{\frac{1}{2}}d(\disp{\max_{0\leq s\leq
t}Y_s})=d(\disp{\max_{0\leq s\leq t}Y_s})$. It follows, from a
result by R.A. Doney, J. Warren and M. Yor \cite{dwy} that this
latter equation has a unique solution in law, and then $X$ is unique
in law.
\medskip

We now consider the following SDE with reflecting boundary and a
measurable drift $b$, i.e.,
\begin{equation}
\label{eqint2} \left\{
\begin{array}{lll}
X_{t} &=&W_{t}+\int_{0}^{t}b(s,X_{s})ds+\alpha \disp{\max_{0\leq
s\leq
t}X_{s}}\,\,+L_{t}^{0},\\
X_{t} &\geq &0,\,\forall t\geq 0.
\end{array}
\right.
\end{equation}
We then have the following result:
\begin{prop}\label{ex}
There exists a weak solution of (\ref{eqint2}) which is moreover
unique in law.
\end{prop}

\noindent {\bf Proof}: Let $(\bar{X}_{t},\bar{W})_{t\geq 0}$ be a
solution in law for (\ref{eqint1}) on a filtered probability space
$(\bar{\Omega}, \bar{\cal F}, (\bar{\cal F}_t)_{t\geq 0},\bar{P})$.
Using It\^o's formula with $(\bar{X})^p$ for $p\geq 2$ we deduce
that for any $t\geq 0$, $\E [\sup_{s\leq t}|\bar{X}_s|^p]<\infty$.
This property combined with the fact $b$ is of linear growth imply
that for any $t\geq 0$ we have :
\begin{equation*}
\bar{\E}\left[\exp \{ -\int_{0}^{t}b(s,
\bar{X}_{s})d\bar{W}_{s}-\int_{0}^{t}b^{2}(s,\bar{X}_{s})ds\}\right]=1.
\end{equation*}
One can see e.g. the appendix of \cite{ElkHam} for the proof. Then
there exists a probability $Q$ on $\bar{\Omega}$ such that the
process
$\tilde{W}:=\{\bar{W}_{t}+\int_{0}^{t}b(s,\bar{X}_{s})ds\}_{t\geq
0}$ is a $Q$-Brownian motion on $\bar{\Omega}$. It follows that the
process $(\bar{X}_{t})_{t\geq 0}$ satisfies the following SDE:
$$
\left\{
\begin{array}{lll}
\bar{X}_{t} &=&\tilde{W}_{t}+\ig{0}{t}b(s,\bar{X}_{s})ds+\al
\max_{0\leq
s\leq t}\bar{X}_{s}\,\,+L_{t}^{0}\\
\bar{X}_{t} &\geq &0,\forall\,\,t\geq 0.
\end{array}\right.
$$
Finally uniqueness in law follows from Girsanov's theorem.
\medskip

Were are now ready to give the main result of this section.
\medskip

\begin{thm}
For any $\alpha \in (0,1) $ the following equation:
\begin{equation}
\label{eqint3} \left\{
\begin{array}{l}
X_{t}=\ig{0}{t}\si (s,X_{s})dW_{s}+\ig{0}{t}b(s,X_{s})ds+\al
\max_{0\leq s\leq
t}X_{s}\,\,\,+L_{t}^{0}\left( X\right),\\
X_{t}\geq 0, \forall t\geq 0.
\end{array}
\right.
\end{equation}
has a unique strong solution which has pathwise uniqueness property
\end{thm}

\noindent \textbf{Proof}: The existence/uniqueness in law is a
consequence of Proposition \ref{ex} and the condition [ND] on $\sigma$,
therefore we only prove pathwise uniqueness.

So suppose that $X$ and $Y$ are two solutions for the SDE
(\ref{eqint3}) defined on a common probability basis with respect to
the same Brownian motion $W$. The assumptions [ND] and [BV2] on
$\sigma $ imply that $L_{.}^{0}( X-Y) =0$ (see e.g. \cite{RY}).
\medskip

We now show that $X\wedge Y$ and $X\vee Y$ are also solutions to
(\ref{eqint3}). By Tanaka's formula we have:
\begin{equation*}
\begin{array}{lll}
X_{t}\vee Y_{t} & = & Y_{t}+( X_{t}-Y_{t}) ^{+}\\
& = & Y_{t}+\int_{0}^{t}1_{\{ X_{s}>Y_{s}\} }d(
X_{s}-Y_{s}) +\frac{1}{2}L_{t}^{0}( X-Y) \\
& = & \int_{0}^{t}1_{\{ X_{s}>Y_{s}\}
}dX_{s}+\int_{0}^{t}1_{\{ X_{s}\leq Y_{s}\} }dY_{s} \\
& = & \int_{0}^{t}\sigma (s,X_{s}\vee
Y_{s})dW_{s}+\int_{0}^{t}b(s,X_{s}\vee Y_{s})ds \\
&  & \qquad +\alpha ( \int_{0}^{t}1_{\{ X_{s}>Y_{s}\}
}dM_{s}^{X}+\int_{0}^{t}1_{\{ X_{s}\leq Y_{s}\}
}dM_{s}^{Y}) \\
&  &  \qquad +\int_{0}^{t}1_{\{ X_{s}>Y_{s}\} }dL_{s}^{0}( X)+
\int_{0}^{t}1_{\{ X_{s}\leq Y_{s}\} }dL_{s}^{0}( Y),
\end{array}
\end{equation*}
where $M_{t}^{Z}$ denotes $\max_{0\leq s\leq t}Z_{s}$. Now from a
result by Ouknine \cite{ouk} (or, Ouknine $\&$  Rutkowski \cite{or}), we
have:
\begin{equation*}
\int_{0}^{t}1_{\left\{ X_{s}>Y_{s}\right\} }dL_{s}^{0}\left(
X\right) +\int_{0}^{t}1_{\left\{ X_{s}\leq Y_{s}\right\}
}dL_{s}^{0}\left( Y\right) =L_{t}^{0}\left( X\vee Y\right).
\end{equation*}
But,
\begin{equation*}
M_{t}^{X\vee Y}=\int_{0}^{t}1_{\left\{ X_{s}>Y_{s}\right\}
}dM_{s}^{X}+\int_{0}^{t}1_{\left\{ X_{s}\leq Y_{s}\right\}
}dM_{s}^{Y}.
\end{equation*}
Consequently, we get that $X\vee Y$ is actually a solution to (\ref{eqint3}).
Similarly we can show that  $X\wedge Y$is a solution as well. As
$X$ and $Y$ have integrable paths on finite time interval, then:
\begin{equation*}
E[| X-Y|] =E[X\vee Y-X\wedge Y]=0
\end{equation*}
which implies that $X=Y$ whence pathwise uniqueness. $\diamond$

\section{SDE involving maximum process in the drift}\label{section3}

In this section we give an existence and pathwise uniqueness result
for SDEs where the drift is a function of the maximum process. Those
equations appear in recent result by B. Roynette P. Vallois and M.
Yor \cite{rvy} see also J. Obloj and M. Yor \cite{oy} about non
Markovian process satisfying the $\mathbf{2M-X}$ property.
\medskip

More precisely, let us consider the following SDE:
\begin{equation}
\label{eqint4} X_{t}=x+W_{t}+\int_{0}^{t}b(X_{s},\max_{0\leq s\leq
t}X_{s})ds,
\end{equation}
where $b$ is measurable and bounded function on $\mathbf{R\times
R\,}$\
to $%
\mathbf{R.}$ We have the following result:

\begin{thm}
Assume that for each $x\in \mathbf{R,}$ $y\rightarrow b\left(
x,y\right) $ is strictly increasing on the set $\left\{ y\in
\mathbf{R/}y\geq x\right\} $. Then (\ref{eqint4}) has a strong
solution which is pathwise unique.
\end{thm}

\noindent {\bf Proof}: The existence and uniqueness in law is a
consequence of Girsanov's theorem. Suppose now that $X$ and $Y$ are
two solutions for SDE (\ref{eqint4})  defined on the same
probability basis with respect to the same Brownian motion $W$, then
$\ X-Y$ is continuous with bounded variation and thus
$L_{.}^{0}\left( X-Y\right) =0.$

Next applying Tanaka's formula yields:
\begin{equation}
\label{eqint5}
\begin{array}{lll}
X_{t}\vee Y_{t} & = & Y_{t}+\left( X_{t}-Y_{t}\right) ^{+}\text{ } \\
& = & Y_{t}+\int_{0}^{t}1_{\left\{ X_{s}>Y_{s}\right\} }d(
X_{s}-Y_{s}) +\frac{1}{2}L_{t}^{0}( X-Y)  \\
& = & x+\int_{0}^{t}1_{\left\{ X_{s}>Y_{s}\right\}
}dX_{s}+\int_{0}^{t}1_{\left\{ X_{s}\leq Y_{s}\right\} }dY_{s} \\
& = & x+W_{t}+\int_{0}^{t}1_{\left\{ X_{s}>Y_{s}\right\}
}b(X_{s},M_{s}^{X})ds \\
&  & \qquad \qquad \qquad +\int_{0}^{t}1_{\left\{ X_{s}\leq
Y_{s}\right\} }b(Y_{s},M_{s}^{Y})ds.
\end{array}
\end{equation}
Observe that by a symmetry argument we have: $\forall t\geq 0$,
\begin{equation*}
\int_{0}^{t}1_{\left\{ X_{s}=Y_{s}\right\}
}b(Y_{s},M_{s}^{Y})ds=\int_{0}^{t}1_{\left\{ X_{s}=Y_{s}\right\}
}b(X_{s},M_{s}^{X})ds.
\end{equation*}
We will now prove that if $X_{s}>Y_{s}$ then $M_{s}^{X}\geq
M_{s}^{Y}.$ Actually let us set
\begin{equation*} s_{0}=\sup \left\{
v\in \left[ 0,s\right] ,\,X_{v}=Y_{v}\right\}.
\end{equation*}
Since $X_{s}-Y_{s}$ is absolutely continuous then
$\dfrac{d}{ds}\left( X_{s}-Y_{s}\right) _{s=s_{0}}\geq 0$ and thus
\begin{equation*}
b(X_{s_{0}},M_{s_{0}}^{X})\geq b(Y_{s_{0}},M_{s_{0}}^{Y}).
\end{equation*} Since $(b(x,y))_{y\geq x}$ is strictly increasing
and through the definition of $s_0$ we deduce that
$M_{s_{0}}^{X}\geq M_{s_{0}}^{Y}$. Furthermore  for $v\in \left[
s_{0},s\right] $ we have $X_{v}> Y_{v}$ and then $M_{s}^{X}\geq
M_{s}^{Y}$ which is the desired result.
\medskip

Going back to (\ref{eqint5}) we obtain:
\begin{equation*}
\begin{array}{lll}
X_{t}\vee Y_{t} & = & x+W_{t}+\int_{0}^{t}1_{\left\{
X_{s}>Y_{s}\right\} }b(X_{s},M_{s}^{X})ds+\int_{0}^{t}1_{\left\{
X_{s}\leq
Y_{s}\right\} }b(Y_{s},M_{s}^{Y})ds \\
& = & x+W_{t}+\int_{0}^{t}1_{\left\{ X_{s}>Y_{s}\right\}
}b(X_{s},M_{s}^{X})ds+\int_{0}^{t}1_{\left\{ X_{s}<Y_{s}\right\}
}b(Y_{s},M_{s}^{Y})ds \\
&  & \qquad\qquad +\int_{0}^{t}1_{\left\{ X_{s}=Y_{s}\right\}
}b(Y_{s},M_{s}^
{Y})ds \\
& = & x+W_{t}+\int_{0}^{t}1_{\left\{ X_{s}>Y_{s}\right\}
}b(X_{s}\vee
Y_{s},M_{s}^{X\vee Y})ds \\
&  & \qquad +\int_{0}^{t}1_{\left\{ X_{s}<Y_{s}\right\} }b(X_{s}\vee
Y_{s},M_{s}^{X\vee Y})ds+\int_{0}^{t}1_{\left\{ X_{s}=Y_{s}\right\}
}b(Y_{s},M_{s}^{Y})ds \\
& = & x+W_{t}+\int_{0}^{t}b(X_{s}\vee Y_{s},M_{s}^{X\vee
Y})ds-\int_{0}^{t}1_{\left\{ X_{s}=Y_{s}\right\}
}b(Y_{s},M_{s}^{X\vee Y})ds
\\
&  & \qquad +\int_{0}^{t}1_{\left\{ X_{s}=Y_{s}\right\}
}b(Y_{s},M_{s}^{Y})ds.
\end{array}
\end{equation*}
But
\begin{equation*}
\begin{array}{ccc}
\int_{0}^{t}1_{\left\{ X_{s}=Y_{s}\right\}
}b(Y_{s},M_{s}^{Y})ds-\int_{0}^{t}1_{\left\{ X_{s}=Y_{s}\right\}
}b(X_{s},M_{s}^{X\vee Y})ds& =
&  \\
\qquad \int_{0}^{t}1_{\left\{ X_{s}=Y_{s}\right\} \cap \left\{
M_{s}^{X}>M_{s}^{Y}\right\} }\left(
b(Y_{s},M_{s}^{Y})-b(X_{s},M_{s}^{X})\right) ds & = & 0.
\end{array}
\end{equation*}
This proves that $X\vee Y$ is also a solution for equation
(\ref{eqint4}). Now uniqueness in law implies that pathwise
uniqueness holds, whence the desired result. $\diamond$
\section{Some possible extensions}\label{section4}
{\bf 5.1}\hskip2mm The perturbed SDE (\ref{eq1}) may be
regarded as a particular case of the the following doubly perturbed
SDE
\begin{eqnarray}\label{e2}
X_t=\xi+\int_0^t \si(s,X_s)dW_s +\int_0^t b(s,X_s)ds +\al\max_{0\leq
s\leq t}X_s +\be \min_{0\leq s \leq t}X_s,
\end{eqnarray}
where $\al,\be \in \R.$ \\
As it was mentioned in the introduction, this doubly perturbed SDE
was the subject of studies of several authors (see the reference in
the beginning of the introduction). Carmona, Petit and Yor (\cite{CPY}) shows, in the particular case when $b\equiv 0$ and $\si$ $\equiv 1$,
the existence and pathwise uniqueness of the strong solution by using a fixed point argument, under the following conditions on the parameters $\al$ and $\be:$
\begin{eqnarray}\label{e3}
\left\{
\begin{array}{ll}
\al < 1,\ \be <1 \\
\disp{\frac{|\al\be|}{(1-\al)(1-\be)}} <1.
\end{array}
\right.
\end{eqnarray}
Here we get strong existence and pathwise uniquness of the solution
to (\ref{e2}), in the Lipschitz case by using the Picard iteration.
More precisely, if $\si$ and $\be$ are taken to satisfy:
\medskip
\begin{eqnarray}\label{e1}
\left\{
\begin{array}{ll}
\left|\si(s,x)-\si(s,y)\right|& \leq c|x-y|\\
\left|b(x)-b(y)\right| & \leq c|x-y|
\end{array}
\right.
\end{eqnarray}
for every $x,y$ in $\R$,\,$s\geq 0$ and some constant $c>0$, then we
have :
\begin{thm}\label{liph}
Under the assumptions (\ref{e3}) and (\ref{e1}), and if the random
variable $\xi$ is such that $\E(|\xi|^2)< \infty $ , then there
exists a unique continuous $\mathcal{F}^{W}$ adapted process $(X_t
,t\geq 0)$ which is a solution to the doubly perturbed SDE (\ref{e2}). Moreover
$\E(\disp\max_{0\leq s\leq T}|X_s|^{2})< \infty$ for every $T>0$.
\end{thm}
\noindent \textbf{Proof}: \emph{Proof of the existence}: The construction of the solution uses the
Picard iteration. Let us consider the sequence implicitly defined
by:
\begin{align}\label{e4}
\left\{
\begin{aligned}
X^0_t & =\frac{\xi}{1-\al },\\
X_t^{n+1}& =\xi + \int_0 ^t \si(X^n_s)dW_s +\int_0^t b(X^n_s)ds +
\al M_t^{n+1} +\beta I^{n+1}_t,
\end{aligned}
\right.
\end{align}
where $M_t^{n+1}\triangleq \disp{\max_{0\leq s\leq t}{X^{n+1}_s}}$
and $I_t^{n+1}\triangleq \disp{\min_{0\leq s\leq t}{X^{n+1}_s}}$ for
every $t\geq 0$.\\
Observe that $(X^n)_{n\geq0}$ is well defined. In fact $X^{n+1}$ is
explicitly evaluated from $X^n$, since by Skorohod's lemma
\cite{RY}, one can easy seen that:
\begin{eqnarray}\label{e5}
\left\{
\begin{array}{ll}
(1-\al)M_t^{n+1} & =\disp{\max_{0\leq s\leq t}{\left(\xi+\ig{0}{s}
\si(X^n_u)dW_u
+\ig{0}{s} b(X^n_u)du + \be I^{n+1}_s\right)}},\\
(\be - 1)I_t^{n+1}  & =\disp{\max_{0\leq s\leq
t}{\left(-\xi-\ig{0}{s} \si(X^n_u)dW_u -\ig{0}{s} b(X^n_u)du - \al
M_s^{n+1}\right)}}.
\end{array}
\right.
\end{eqnarray}
Combining these two equalities, we get:
\begin{eqnarray}\label{e6}
M_t^{n+1}&=&\frac{1}{1-\al}\max_{0\leq s\leq
t}\left\{\left(\xi+\ig{0}{s}
\si(X^n_u)dW_u +\ig{0}{s} b(X^n_u)du \right)\right.\nonumber\\
&+&\left.\frac{\be}{\be-1}\max_{0\leq u\leq s}\left(-\xi-\int_0 ^u
\si(X^n_v)dW_v - \ig{0}{u} b(X^n_v)dv - \al M_u^{n+1}
\right)\right\}.\nonumber\\
\end{eqnarray}
By a similar argument as in \cite{CPY}, using a fixed
point theorem and the hypothesis (\ref{e3}), we get the well
definiteness and the adaptation of $M_t^{n+1}$ with respect to the
filtration of $B$. Thus, $(X^n)_{n\geq0}$
is well defined.\\
Now, let us show that $X^n$ converges uniformly on compact
intervals almost surely. On one hand, we have:
\begin{eqnarray}\label{e6}
\max_{0\leq s\leq t}|X^{n+1}_s -X^n_s|\leq a_n (t)+ b_n (t)+
|\al|\max_{0\leq s\leq t}|M^{n+1}_s -M^n_s| + |\be|\max_{0\leq s\leq
t}|I^{n+1}_s -I^n_s|.
\end{eqnarray}
On the other hand, by (\ref{e5})
\begin{align}\label{e8}
\left\{
\begin{aligned}
\max_{0\leq s\leq t}|M^{n+1}_s -M^n_s| & \leq\frac{1}{1-\al }\{a_n
(t)+ b_n (t)\}+\frac{|\be|}{1-\al}\max_{0\leq s\leq t}|I^{n+1}_s
-I^n_s|,\\
\max_{0\leq s\leq t}|I^{n+1}_s -I^n_s| & \leq\frac{1}{1-\be }\{a_n
(t)+ b_n (t)\}+\frac{|\al|}{1-\be}\max_{0\leq s\leq t}|M^{n+1}_s
-M^n_s|,\\
\end{aligned}
\right.
\end{align}
where $a_n(t):=\disp{\max_{0\leq u\leq
t}|\int_0^u(\si(X^n_s)-\si(X^{n-1}_s))dW_s)|}$
and $b_n (t):=\disp{\int_0^t |b(X^n_s)- b(X^{n-1}_s)|ds}$,\, for $t\geq 0$. Thus,
\begin{align}\label{e7}
\left(1-\frac{|\al\be|}{(1-\al)(1-\be)}\right)\max_{0\leq s\leq
t}|M^{n+1}_s -M^n_s|  \leq
\frac{1}{1-\al}\left(1+\frac{|\be|}{1-\be}\right)(a_n (t)+ b_n (t)).
\end{align}
 Combining (\ref{e6}) and (\ref{e7}) yields
\begin{eqnarray}
\max_{0\leq s\leq t}|X^{n+1}_s -X^n_s|& \leq &
\left(1+\frac{|\be|}{1-\be}\right)\left(1+\frac{|\al|
(1-\be)}{(1-\al)(1-\be)-|\al\be|}\right.\nonumber\\
&+&\left.\frac{|\al\be|}{(1-\al)(1-\be)-|\al\be|}\right)(a_n (t)+b_n
(t)).\nonumber
\end{eqnarray}
Applying Bukholder-Davis-Gundy's inequality and Lipshitz's
conditions to this inequality, we obtain for a generic constant
$c>0$,
\begin{align}
\disp{\E \left(\max_{0\leq s\leq t}{|X^{n+1}_s -X^n_s|^2}\right)\leq
c\int_0^t\E|X^{n}_s -X^{n-1}_s|^2 ds}.
\end{align}
 Iterating this inequality led, for every $T>0$, to
\begin{eqnarray}
\E \left[\max_{0\leq s\leq t}{|X^{n+1}_s -X^n_s|^2}\right]\leq
c\frac{T^n}{n!}
\end{eqnarray}
since\quad $\E\xi^2<\infty$. Consequently, by the Chebychev's
inequality and the Borel-Cantelli lemma, we get the uniform
convergence of the sequence $(X^{n})_n$ to a continuous process $X$
on $[0,T]$(see for instant \cite{SC}). Letting $n\rightarrow\infty$
in (\ref{e4}) it follow that $X$ is a solution to the doubly
perturbed SDE (\ref{e2}). As $T$ is arbitrary
this proves the existence.\\
\emph{Proof of the uniqueness}: Now suppose that $X$ and $Y$ are two solution to the SDE
 (\ref{e2}) with some initial condition $\xi$ and some driving
 Brownian motion $B$, then we can easily seen
 \begin{equation}\label{e9}
    |X_t - Y_t|\leq |\int_0^t(\si(X_s)-\si(Y_s))dW_s| + |\be||I^X_t
    -I^Y_t| +\frac{|\al|}{1-\al}\max_{0\leq s\leq t}|(A^{X}_{s})^{+}- (A^{Y}_{s})^{+}|
 \end{equation}
 with  $A^{\om}_{t}:= \xi+\int_0^t \si(s,\om_s)dB_s +\int_0^t b(s,\om_s)ds +\al
 I^X_t$ where $\om$ is $X$ or $Y$.\\
 Arguing as above (using Skorohod's lemma), we have
 \begin{align}\label{e10}
\left(1-\frac{|\al\be|}{(1-\al)(1-\be)}\right)\max_{0\leq s\leq
t}|I^{X}_s -I^Y_s|  \leq c(a (t)+ b (t)),
\end{align}
where $a(t):=\disp{\max_{0\leq u\leq
t}|\int_0^u(\si(X_s)-\si(Y_s))dW_s)|}$
  and $b(t):=\disp{\int_0^t |b(X_s)- b(Y_s)|ds}$.\\
 Moreover,
 \begin{equation*}
   |A^{X}_{t}-A^{Y}_{t} |\leq|\ig{0}{t}(\si(X_s)-\si(Y_s))dW_s|+ \ig{0}{t} |b(X_s)-
   b(Y_s)|ds +|\be||I^{X}_t -I^Y_t|
 \end{equation*}
which implies, taking into account of (\ref{e10})
\begin{equation}\label{e11}
  \max_{0\leq s\leq t}|A^{X}_{s}- A^{Y}_{s}|\leq c(a (t)+ b (t))
\end{equation}
and consequently
\begin{equation}\label{e12}
    |X_t - Y_t| \leq c(a (t)+ b (t)).
\end{equation}
Finally, by applying BDG's inequality and the lipshitz condition we
get $$\E (|X_t - Y_t|^{2})\leq c\ig{0}{t}\E\left(|X_u -
Y_u|^{2}|\right)du.$$
Hence, $\E\left(|X_t-Y_t|^{2}\right)=0$ by Gronwall's lemma. Thus the
solution is unique. $\diamond$\\
\begin{remark} It would be interesting to see, if the result of Theorem.\ref{liph} can be proved by the approach used in the first section of this paper. This will illuminate the situation when the coefficients $\si$ and $b$ are not
Lipschitz.
\end{remark}

\noindent
{\bf 5.2}\hskip2mmLet $H$ be an absolutely continuous and
increasing function such that $0<H^{\prime }\left( x\right) <1$ and
$H(0)=0$. We can develop the same arguments as previously to show
that our results on strong existence and pathwise uniqueness are still
valid for the following SDEs:
\begin{equation}
X_{t}=\displaystyle{\ig{0}{t}\si(s,X_{s})dW_{s}+\ig{0}{t}b(s,X_{s})ds+H\left(
\max_{0\leq s\leq t}X_{s}\right)}, \,\forall t\geq 0
\end{equation}or
\begin{equation}
\left\{
\begin{array}{lll}
X_{t} &=&X_{0}+\ig{0}{t}\si
(s,X_{s})dW_{s}+\ig{0}{t}b(s,X_{s})ds+H\left( \disp{\max_{0\leq
s\leq
t}X_{s}}\,\right) \,\,+L_{t}^{0} \\
X_{t} &\geq &0,\,\,\forall t\geq 0
\end{array}
\right.
\end{equation}
and finally
\begin{equation}
X_{t}=X_{0}+\ig{0}{t}\sigma
(s,X_{s})dW_{s}+\ig{0}{t}b(s,X_{s})ds+\ig{0}{t}\al( s)
dM_{s}^{X}\,\,\forall t\geq 0,
\end{equation}
for any deterministic Borel function $\alpha $ such that \ $0$ $\leq
\alpha \left( s\right) <1.$ $\diamond$
\medskip


\end{document}